\documentclass{cmslatex}

\usepackage{latexsym, amssymb, enumerate, amsmath}
\usepackage{amssymb}
\usepackage{graphicx}
\usepackage{color}
\usepackage[update,prepend]{epstopdf}

\renewcommand{\theequation}{\arabic{section}.\arabic{equation}}

\newcommand{\bE}{{\mathbb E}}
\newcommand{\bP}{{\mathbb P}}
\newcommand{\eps}{{\varepsilon}}

\title{Stochastic Mode-Reduction in Models with Conservative Fast Sub-Systems}
\author{
Ankita Jain\thanks{
University of Notre Dame,
Department of Applied and Computational Mathematics and Statistics,
153 Hurley Hall, Notre Dame, IN 46556, 
Ankita.Jain@nd.edu}
\and
Ilya Timofeyev\thanks{
University of Houston,
Department of Mathematics,
651 Philip G. Hoffman Hall, 
Houston, TX 77204-3008, 
ilya@math.uh.edu}, 
\and 
Eric Vanden-Eijnden\thanks{
New York University,
Courant Institute,
251 Mercer Street,
New York, NY 10012, 
eve2@cims.nyu.edu}}

\begin{document}
\maketitle


\begin{abstract}
  A stochastic mode reduction strategy is applied to multiscale
    models with a deterministic energy-conserving fast sub-system.
  Specifically, we consider situations where the slow variables are
driven stochastically and interact with the fast sub-system in an
energy-conserving fashion.  Since the stochastic terms only affect the
slow variables, the fast-subsystem evolves deterministically on a
sphere of constant energy. However, in the full model the radius of
the sphere slowly changes due to the coupling between the slow and
fast dynamics. Therefore, the energy of the fast sub-system
  becomes an additional hidden slow variable that must be accounted
for in order to apply the stochastic mode reduction technique to
systems of this type.
\end{abstract}

\section{Introduction}

The derivation of reduced models for various complex high-dimensional
systems has been an active area of research for many decades, and it
has recently received increased attention with particular emphasis on
practical applications. For example, a stochastic mode reduction
strategy was proposed in the context of atmospheric applications in
\cite{mtv2,mtv4}. This strategy has been successfully applied to
several prototype problems \cite{mtv3,mtv5,bati09}, and it also has
been applied to more realistic geophysical examples
\cite{frmava05,frma06}.  The technique builds on the idea of adiabatic
elimination of fast variables in multiscale stochastic systems
\cite{kha66a,kha66b,kur73,pap76} and assumes the existence of scale
separation between the (low-dimensional) slow and (high-dimensional)
fast variables. In the context of atmospheric applications, it is
plausible to consider systems where the fast variables are driven
stochastically and interact with the slow dynamics in an
energy-conserving fashion.  Since then, this technique has been
extended to purely deterministic conservative systems with a chaotic
deterministic heat bath constituting the fast variables
\cite{mtv6,melst11}.

In this paper we present a mathematical formalism which extends the
stochastic mode reduction to systems where the fast sub-system is
deterministic and energy-conservative, but the slow variables can
evolve in a stochastic or chaotic fashion. Therefore, unlike in
\cite{mtv6}, the total energy of the system is not conserved and the
interactions of the slow variables with the fast bath affect the total
energy of the fast sub-system.  The most striking examples of such
models are coupled parabolic-hyperbolic systems where the hyperbolic
part plays the role of the fast sub-system. Such models arise, for
example, in thermo-elasticity \cite{pderef-book,pderef4,pderef5} and
thermo-visco-elasticity \cite{pderef1,pderef2,pderef3}, where they can
be viewed as simplified versions of the true fluid-structure
interaction models: the hyperbolic and the parabolic parts describe
the fast-moving waves and the temperature in the domain, respectively.

Applications of the stochastic mode reduction strategy typically
require understanding the statistical behavior of the fast
sub-system. Often, a mixed analytical-computational approach can be
employed to evaluate the statistical behavior of the fast conservative
sub-system on a shell of constant energy.  However, since the energy
of the fast sub-system is changing in time through the interactions
with the slow modes, it is necessary to explicitly introduce the
energy of the fast sub-system as an additional hidden slow variable in
order to extend the stochastic mode reduction strategy to models
described in the previous paragraph.  
We demonstrate that a closed-form stochastic differential equation for
the evolution of the slow dynamics (including energy) can be derived
and the coefficients in this equation can be evaluated from a
  single numerical simulation of the fast-sub-system on an (arbitrary)
  shell of constant energy.  Here we illustrate our approach on a
  simple prototype model -- more realistic examples will be considered
  in subsequent papers. In particular, we consider an extended triad
model which elucidates analytical and numerical issues of the
stochastic-mode reduction for systems with deterministic conservative
heat bath.

\section{Prototype Triad Model}

To illustrate our approach we consider a prototype extended triad model
\begin{equation}
\begin{aligned}
\dot{x} &=  \frac{1}{\eps} \sum_{j,k=1}^n  A^{xyy}_{1jk} y_j y_k - \gamma x + \sigma \dot{W}, \\
\dot{y_j} &=  \frac{1}{\eps} \sum_{k=1}^n  A^{yxy}_{j1k} x y_k +   \frac{1}{\eps^2} \sum_{i,k=1}^n B^{yyy}_{jik} y_i y_k,
\label{triad}
\end{aligned}
\end{equation}
where $x$ is the slow variable, $y_j$, $j=1,\ldots,n$ are the fast variables, 
$n$ is the total number of fast variables, and
$A^{xyy}_{1jk}$, $A^{yxy}_{j1k}$, $B^{yyy}_{jik}$ are interaction coefficients obeying the
energy-conserving relationships 
\begin{equation}
A^{xyy}_{1jk} + A^{yxy}_{j1k} + A^{yxy}_{k1j} = 0,
\label{c1}
\end{equation}
\begin{equation}
B^{yyy}_{jik} + B^{yyy}_{ikj} + B^{yyy}_{kji} = 0.
\label{c2}
\end{equation}
The fast variables are determined by the $\eps^{2}$ timescale in the
equation \eqref{triad} and the virtual fast sub-system becomes
\begin{equation}
\label{fastsubs}
\dot{y_j} =  \sum_{i,k=1}^n B^{yyy}_{jik} y_i y_k.
\end{equation}
Condition~\eqref{c2} is essential for the formalism presented
  below. It ensures that the virtual fast sub-system in
  \eqref{fastsubs} is energy conservative which guarnatees the
  existence of an invariant measure for this sub-system. The existence
  and uniqueness of this invariant measure is a necessary condition
  for the applicability of the approach presented in this
  paper. Condition \eqref{c1}, on the other hand, ensures that the
  energy is conserved by the interactions between the slow variable
  and the fast modes and, thus, \eqref{c1} guarantees the existence of
  the invariant measure for the full model in \eqref{triad}. We stress
  that we made this assumption for simplicity: the existence of an
  invariant mesaure for the full system~\eqref{triad} is not necessary
  for our approach, and the formalism presented here is also
  potentially applicable to systems where the interaction between the
  slow variables and the fast sub-system is not energy-conservative.

One can attempt to apply the stochastic mode reduction strategy to the
model in \eqref{triad} directly. In the stochastic mode reduction
procedure one needs to consider the behavior of the fast-subsystem
\eqref{fastsubs} and certain statistical properties of $y$-variables
(e.g. correlations functions) typically enter as coefficients into the
reduced model. The behavior of the fast-subsystem depends drastically
on the energy level, but the energy in the fast modes changes in time
due to the interaction with the slow variables. Therefore, it is
necessary to consider the energy of the fast sub-system
\begin{equation}
\label{E}
E = \sum\limits_{k=1}^n y_k^2
\end{equation}
as an additional slow variable to explicitly keep track of the changes
in the statistical behavior of $y$-variables.  The extended multiscale
triad model becomes
\begin{equation}
\label{fullmodel}
\begin{aligned} 
  \dot{x} &= \frac{1}{\eps} \sum_{j,k=1}^n A^{xyy}_{1jk} y_j y_k 
  - \gamma x + \sigma \dot{W},\\
  \dot{E} &= -2 \frac{1}{\eps} \sum_{j,k=1}^n A^{xyy}_{1jk} x y_j y_k,\\ 
\dot{y_j} &= \frac{1}{\eps} \sum_{k=1}^n A^{yxy}_{j1k} x y_k + 
\frac{1}{\eps^2}  \sum_{i,k=1}^n B^{yyy}_{jik} y_i y_k,
\end{aligned}
\end{equation}
where we used the property \eqref{c1} to simplify the equation for the energy.

\subsection{Stationary Distribution of the Triad Model}

The stationary distribution of the generalized triad model
\eqref{fullmodel} can be computed explicitly for any $\eps$ and it is
easy to show that the stationary distribution for $x,y_1,\ldots,y_n$
is a product of Gaussian distributions with mean zero and identical
variances $\sigma^2/(2\gamma)$. In particular, the stationary
distribution does not depend on the parameter $\eps$.

The Fokker-Planck equation for the invariant density $\rho = \rho(x,y_1,..,y_n)$ for 
$x,y_1,\ldots,y_n$ in \eqref{fullmodel} can be written as
\[
0 = \left(A_0 + \frac{1}{\eps}A_1 + \frac{1}{\eps^2} A_2 \right) \rho
\]
where 
\begin{eqnarray} 
A_0 &=& \gamma \partial_{x}x + \frac{\sigma^2}{2} \partial_{xx},
\nonumber \\
A_1 &=& -\sum_{j,k} A^{xyy}_{1jk} y_j y_k \partial_x  - \sum_{j,k} A^{yxy}_{j1k} x \partial_{y_j} y_k,
\nonumber \\ 
A_2 &=& -\sum_{i,j,k} B^{yyy}_{jik} \partial_{y_j} y_i y_k.
\nonumber
\end{eqnarray}

The differential operator $A_0$ annihilates the Gaussian density
\begin{equation}
\label{rhox}
\rho_x(x) = \frac{\sqrt{2 \gamma}}{\sqrt{2 \pi} \sigma} \exp \left(-\frac{\gamma}{\sigma^2} x^2 \right).
\end{equation}
Operators $A_1$ and $A_2$ annihilate separately any function of the
full energy $x^2 + E$ with $E$ in \eqref{E} due to the conservation of
energy properties in \eqref{c1} and \eqref{c2}. Therefore, they also
annihilate the function
\[
\rho = \left( \frac{\sqrt{2 \gamma}}{\sqrt{2 \pi} \sigma} \right)^{n+1}
\exp \left(-\frac{\gamma}{\sigma^2} \left[x^2 + \sum_j y_j^2 \right] \right).
\]

Thus, the invariant measure for $x,\, y_1,\ldots,y_n$ in
\eqref{fullmodel} is a product measure with
$$ 
x, y_1,\ldots,y_n \sim N\left(0,\frac{\sigma^2}{ 2 \gamma}\right).
$$
The stationary distribution for the energy $E$ in \eqref{E} can be
easily derived since the joint stationary distribution of fast
variables is a product of Gaussian densities and, therefore, the
stationary distribution for $E$ is a $\chi$-squared distribution
\begin{equation}
\label{rhoE}
\rho_E(s) = 
C s^{(n-2)/2} e^{- \frac{s\gamma}{\sigma^2}} .
\end{equation}
The joint stationary distribution for the slow variables $x$ and $E$
is a product density of the corresponding marginal densities in
\eqref{rhox} and \eqref{rhoE}.

\subsection{Limit of Full Model as $\eps \rightarrow 0$}

Derivation of the reduced model for $x$ in $E$ in \eqref{fullmodel}
proceeds in a typical fashion.  The Kolmogorov backward equation
associated with \eqref{fullmodel} for a scalar function
$u=u(t,x,E,y_1,\ldots,y_n)$ is given by
\begin{equation} \partial_t u = L_0 u + \frac{1}{\epsilon} L_1
  u + \frac{1}{\epsilon^2} L_2 u,
\label{back}
\end{equation}
where the operators above are  
\begin{eqnarray} 
L_0 &=& - \gamma x \partial_{x} + \frac{\sigma^2}{2} \partial_{xx} ,
\nonumber \\ 
L_1 &=& \sum_{j,k} A^{xyy}_{1jk} y_j y_k \partial_x  + \sum_{j,k} A^{yxy}_{j1k} x y_k \partial_{y_j}  
    - 2 \sum_{j,k} A^{xyy}_{1jk} x y_j y_k \partial_E ,
\nonumber \\ 
L_2 &=& \sum_{i,j,k} B^{yyy}_{jik} y_i y_k  \partial_{y_j} .
\nonumber
\end{eqnarray}

Next, we introduce the projection operator
\[
\bP \, \cdot = \int \, \cdot \,  d \mu(\vec{y}|E), 
\]
where $d\mu(\vec{y}|E)$ is the invariant measure of the fast
sub-system \eqref{fastsubs}.  The energy of the fast sub-system does
not change in time and, thus, the measure $d\mu(\vec{y}|E)$ is
concentrated on the sphere of constant energy. In the absence of
additional information, we assume that the measure $d\mu(\vec{y}|E)$
is the uniform measure on the sphere, i.e.
\begin{equation}
\label{mu}
d\mu(\vec{y}|E) =  S_n^{-1} E ^{1-n/2} \delta(E - \sum y_j^2) \, d\vec{y}, 
\end{equation}
where $n$ is the total number of fast variables, $E$ is the fixed
energy in the fast subsystem, $S_n^{-1} E ^{1-n/2}$ is the normalizing
factor such that $S_n$ does not depend on $E$, and $\delta(\cdot)$ is
the Dirac delta function.

Considering the expansion 
\[
u = u_0 + \eps u_1 + \eps^2 u_2
\]
and collecting powers of $\eps$ we obtain 
\begin{eqnarray}
O(\epsilon^{-2})&:& L_2 u_0 = 0 ,
\nonumber \\ 
O(\epsilon^{-1})&:& L_2 u_1 = - L_1 u_0,
\label{L operators} \\
O(1)&:&  L_2 u_2=  \partial_t u_0 - L_0 u_0 - L_1 u_1.
\nonumber
\end{eqnarray}
The first equation implies that $u_0 = u_0(x,E,t)$ and is independent
of the fast variables $y_1,..,y_n$ since $L_2$ involves
differentiation with respect to only $y$ variables. Applying $\bP$ to
the second equation in \eqref{L operators} we obtain the compatibility
condition $\bP L_1 = 0$ (since $\bP$ is the projection with respect to
the invariant measure of the fast subsystem and, thus $\bP L_2 = L_2
\bP = 0$).  The second equation in (\ref{L operators}) implies that
\begin{eqnarray}
u_1 = - L_2^{-1} L_1 u_0.
\label{u_1 equation}
\end{eqnarray}

Finally, substituting \eqref{u_1 equation} into the third equation and applying  $\bP$ to both sides we obtain the 
reduced (homogenized) equation for $u_0$
\begin{equation}
\label{back reduced}
\partial_t u_0 =L_0 u_0 + L u_0,
\end{equation}
where $L$ is the reduced operator given by
\begin{equation}
\label{L}
L = -\bP L_1 L_2^{-1} L_1\bP.
\end{equation}
Note, that the stochastic terms in the equation for $x$ which are not affected by 
the fast sub-system are simply ``carried through''
by the mode reduction procedure, i.e. the reduced backward equation involves the term $L_0 u_0$,
where operator $L_0$ corresponds to the terms without $\eps$ in the equation for $x$. 
Derivation of the reduced operator proceeds in a formal fashion by manipulating the 
summations in the expression \eqref{L} and formally expressing the action of $L_2^{-1}$. 
Details of the derivation are presented in appendix \ref{a1}.

The compatibility condition $\bP L_1 = 0$ implies the following conditions on the stationary 
first and second moments of the $y$-variables in the statistical behavior of the fast sub-system
 \eqref{fastsubs}
\begin{equation}
\label{ccond}
\bE_\mu y_j = 0, \qquad
\bE_\mu y_j y_k  = 0 {\rm ~for~} j \ne k.
\end{equation}    
The first condition is quite plausible, it simply implies the symmetry of the 
distribution. The second condition is not obvious. While the the stationary second mixed moment $\bE y_j y_k$
is zero in the full model, the requirement in \eqref{ccond} is for the dynamics of the deterministic 
fast sub-system. Therefore, the conditions in \eqref{ccond} need to be verified numerically 
in the simulations of the fast sub-system on a shell of constant energy.

The backward reduced operator can be computed explicitly and after substituting the 
action of $L_2^{-1}$ and some formal manipulations it becomes  
\begin{eqnarray} \nonumber
L &=& - \frac{E^{1/2}}{n^{3/2}} M  \left( (n+1) x \partial_{x} + 2 E \partial_E - 2 x^2 (n+1) \partial_{E} \right) 
\\
& & + \frac{E^{3/2}}{n^{3/2}} M \left(\partial_{xx} - 2 x \partial_{xE} - 2 x \partial_{Ex} + 4 x^2 \partial_{EE} \right),
\label{L final}
\end{eqnarray}
where $n$ is the number of fast variables and $M$ is a numerical constant describing some bulk 
statistical properties of the fast sub-system.
In particular, $M$ is the area under the fourth-order two-point moment of the fast sub-system \eqref{fastsubs}
computed on the energy level $E=n$
\begin{equation}
M = \sum\limits_{jkj'k'} A^{xyy}_{1jk} A^{xyy}_{1j'k'} 
\int\limits_0^\infty \int y_j y_k  y_{j'}(\tau) y_{k'}(\tau) 
\, d\mu(\vec{y}|E=n) \, d\tau .
\label{M}
\end{equation}
The numerical constant $M$ cannot be deduced analytically and has to be computed numerically.
Nevertheless, this computation needs to be done only once using time-averaging in the simulations of the fast sub-system.
Moreover, the fast sub-system alone is not a multiscale system and, thus, can be simulated rather easily.

The effective SDE corresponding to the reduced operator \eqref{L final} is
\begin{eqnarray}
dx &=& - \gamma x dt - (n+1) M x \frac{E^{1/2}}{n^{3/2}} dt + \sigma dW_1 + \sqrt{2 M} \left(\frac{E}{n}\right)^{3/4} dW_2 ,
\nonumber \\
\label{reduced} \\
dE &=&  -2 M \left(\frac{E}{n}\right)^{3/2} dt + 2 (n+1) M x^2 \frac{E^{1/2}}{n^{3/2}} dt - 2\sqrt{2 M}x \left(\frac{E}{n}\right)^{3/4} dW_2 .
\nonumber
\end{eqnarray}
%


\section{Numerical Simulations}

To illustrate the mode reduction approach we choose parameters of the full model
 (\ref{triad}) as
\begin{eqnarray}
\gamma = 1, \quad \sigma = 2.236, \quad  n = {\rm dim}(\vec{y}) = 10.
\label{parameters}
\end{eqnarray}
The interaction coefficients are given in Tables \ref{tabA} and \ref{tabB} in appendix \ref{a2}. 
The corresponding number of triads in the full model is
\[
{\rm Number~of~} xyy {\rm ~triads~} = 10, \qquad
{\rm Number~of~} yyy {\rm ~triads~} = 19.
\]
%

%
We use an analog of the split-step method to integrate the full model
\eqref{triad}.  In particular, due to the multiscale nature of the
full model, the time-stepping for the deterministic part of the full
model requires a relatively high order discretization formula in order
to correctly represent the energy transfer between $x$ and
$y$-variables during a single time-step. We observed empirically that
the full model is rather sensitive to the deterministic integrator;
thus, we use a RK5 formula to integrate numerically the deterministic
part of the system \eqref{triad}. Next, we use Euler discretization to
add a Gaussian random variable which approximates $\Delta W$. The
overall (stochastic) scheme is of order $\Delta t$, but the energy
transfer in this numerical scheme is represented much more accurately.
This turns out to be essential. Indeed we observed empirically
  that numerical simulations of the full system with the low-order
  deterministic discretization (Euler and RK2) produce severe
  discrepancies with analytical predictions (e.g. the stationary
  distribution of the energy variable given by \eqref{rhoE}). Both RK4
  and RK5 discretizations of the deterministic terms produce similar
  results in this case, and we use a higher-order method to ensure the
  robustness of numerical results.

We compute the statistical properties of the slow variable, $x$, as
well as the statistics of $E$, since $E$ enters explicitly into the
reduced model and stationary statistics of both variables should be
used as a criteria for the comparison between the full model and the
reduced dynamics.

We use time-averaging combined with Monte-Carlo approach to accelerate
computations of the full model in \eqref{fullmodel}.  In particular,
we perform $10$ independent runs with different initial conditions of
the full model with $T = 40,000$ for time-averaging in each run and
then average over these $10$ realizations.  We use the time-step
$\Delta t = 10^{-4}, \, 2.5 \times 10^{-5}, \, 2\times 10^{-5}, \,
10^{-6}\, $ for $\eps = 1, \, 0.5, \, 0.25, \, 0.1$, respectively.

The reduced model is integrated using the same method with the RK5 deterministic integrator and 
Euler discretization for the noise. 
Compared to the full model, the reduced system is not as sensitive to the choice of the discrete 
scheme for the deterministic terms. In particular, the numerical results with the RK2 and RK5 discretizations 
are very similar. However, a  much smaller time-step is required for 
the numerical scheme with the Euler discretization (i.e. the stochastic Euler-Maruyama scheme) to 
produce accurate results. 
We also perform a hybrid approach combining time-averaging and Monte-Carlo simulations
to compute the stationary statistical properties of the reduced model. 
We run $K=10$ trajectories each with $T = 100,000$,  $\Delta t = 10^{-5}$, and different 
initial conditions, and then take the average of statistical quantities. 
The coefficient $M$ in the reduced model \eqref{reduced} computed from the fast subsystem 
numerically is
$$
M = 1.2759.
$$

\begin{figure}[ht]
\centerline{\includegraphics[width=12cm]{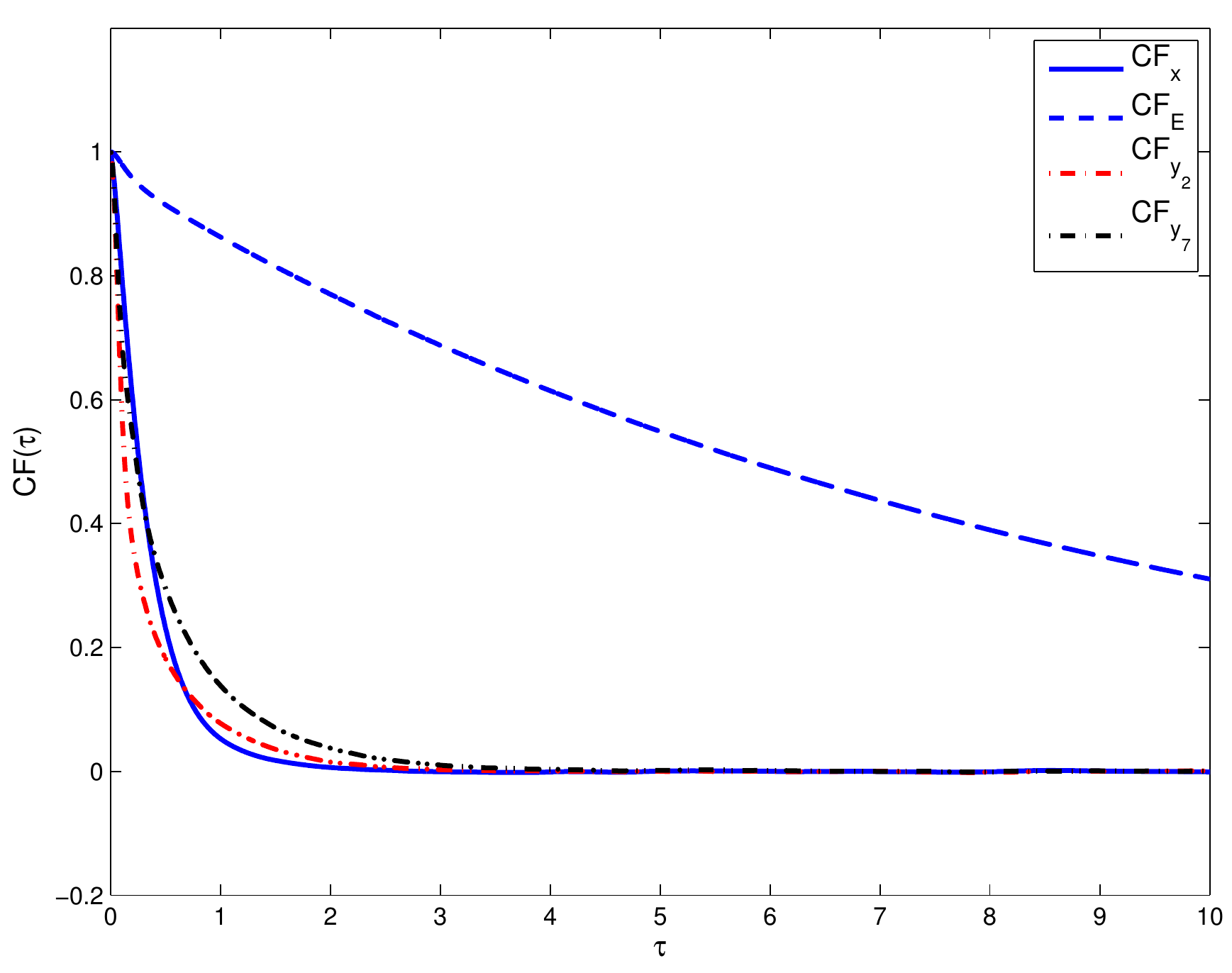}} 
\caption{Normalized correlation functions  $CF_x$, $CF_E$ and $CF_{y_2}$ and $CF_{y_7}$; 
modes $y_2$ and $y_7$ are among the slowest $y$-variables and their correlation functions 
illustrate the absence of scale separation for $\eps=1$.}
\label{fig1}
\end{figure}

First, we examine the scale separation between the slow and the fast variables in the full model \eqref{triad}. 
Comparison of correlation functions of the slow variables and a few fast modes for $\eps=1$ is presented in Figure \ref{fig1}. 
The corresponding correlation times for $\eps=1$ are presented in Table \ref{ct}. 
Correlation times are computed as the inverse of the area under the graph of the normalized 
correlation function.
This simulation shows that there is a group of $y$-variables with time-scales comparable to the
time-scale of $x$. Therefore, we can conclude that there is no time-scale separation between 
the $x$ and $y$-variables for $\eps=1$. 
However, we also observe that $E$ is much slower than $x$.
Moreover, Table \ref{ct} demostrates that the scale separation between the slow variables $(x,E)$
and fast $y$-variables increases as $\eps\to 0$.
Therefore, although we expect that the stochastic mode reduction will perform quite well in the limit $\eps\to0$,
the agreement with the numerical results for $\eps=1$ is not guaranteed.
In addition, Figure \ref{fig1} also demonstrates the necessity to include the energy $E$ in the reduced description.
The energy $E \equiv E(t)$ acts as a hidden stochastic variable which is essential for the derivation of the reduced 
model; 
Figure \ref{fig1} indicates that it is not possible, for instance, to average the reduced model \eqref{reduced} 
with respect to the stationary distribution for $E$, since $E$ is much slower than $x$.
\begin{table}
\label{ct}
\begin{center}
    \begin{tabular}{ | l | c |c|c|c|}
    \hline
    Variable & CT $\eps=1$ & CT $\eps=0.5$ & CT $\eps=0.25$ & CT $\eps=0.1$ \\ \hline
    $x$ & 0.3395 & 0.3433 & 0.3324 & 0.3258 \\ \hline
    $E$ & 8.0605 &  7.7616 & 7.336  & 7.218   \\ \hline
    $y_1$ &  0.2395 &  0.0645  & 0.0177 & 0.0052 \\ \hline
    $y_2$ &  0.3239 &  0.0828  & 0.022   & 0.006   \\ \hline
    $y_3$ &  0.3148 &  0.088    & 0.025   & 0.0063  \\ \hline
    $y_4$ &  0.1460 & 0.0403   & 0.0125  & 0.0045  \\ \hline
    $y_5$ &  0.1706 & 0.0461   & 0.0131  & 0.0046  \\ \hline
    $y_6$ &  0.0713 & 0.0209   & 0.0068  & 0.0041  \\ \hline
    $y_7$ &  0.4910 & 0.1255   & 0.0039  & 0.0076  \\ \hline
    $y_8$ &  0.0570 & 0.0156   & 0.0056  & 0.0041 \\ \hline
    $y_9$ &  0.3383 & 0.0876   & 0.0244  & 0.006  \\ \hline
    $y_{10}$ & 0.1126 & 0.0305 & 0.0094  & 0.0041 \\
    \hline
    \end{tabular}
\end{center}
\caption{Correlation Time (CT) in the simulations of the full model with different values of $\eps$.}
\end{table}

\begin{figure}[ht]
\centerline{\includegraphics[width=12cm]{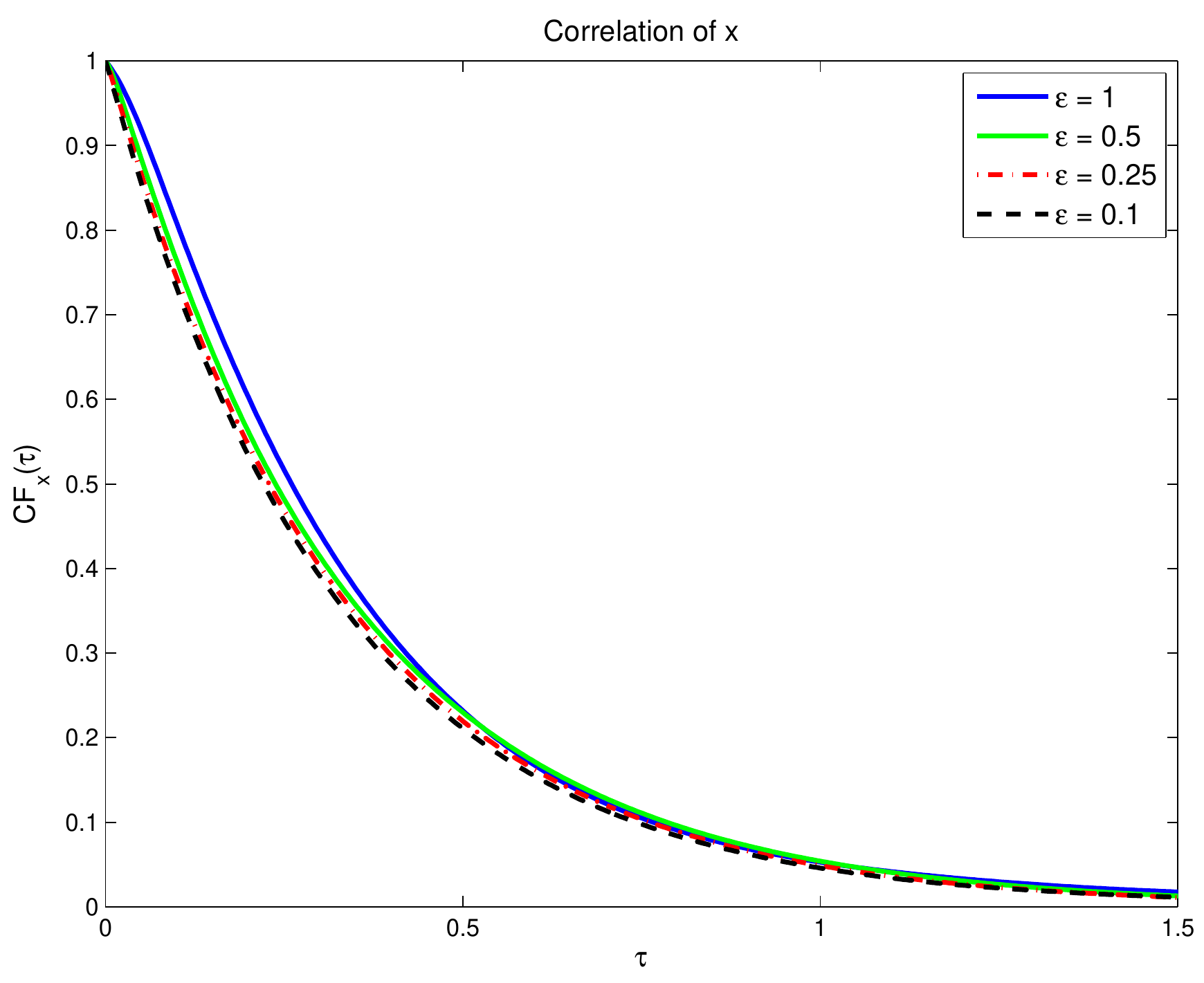}} 
\caption{Normalized correlation function of $x$ in the full model \eqref{triad} 
for $\eps=1, \, 0.5, \, 0.25, \, 0.1$.}
\label{fig:cfxfull}
\end{figure}

\begin{figure}[ht]
\centerline{\includegraphics[width=12cm]{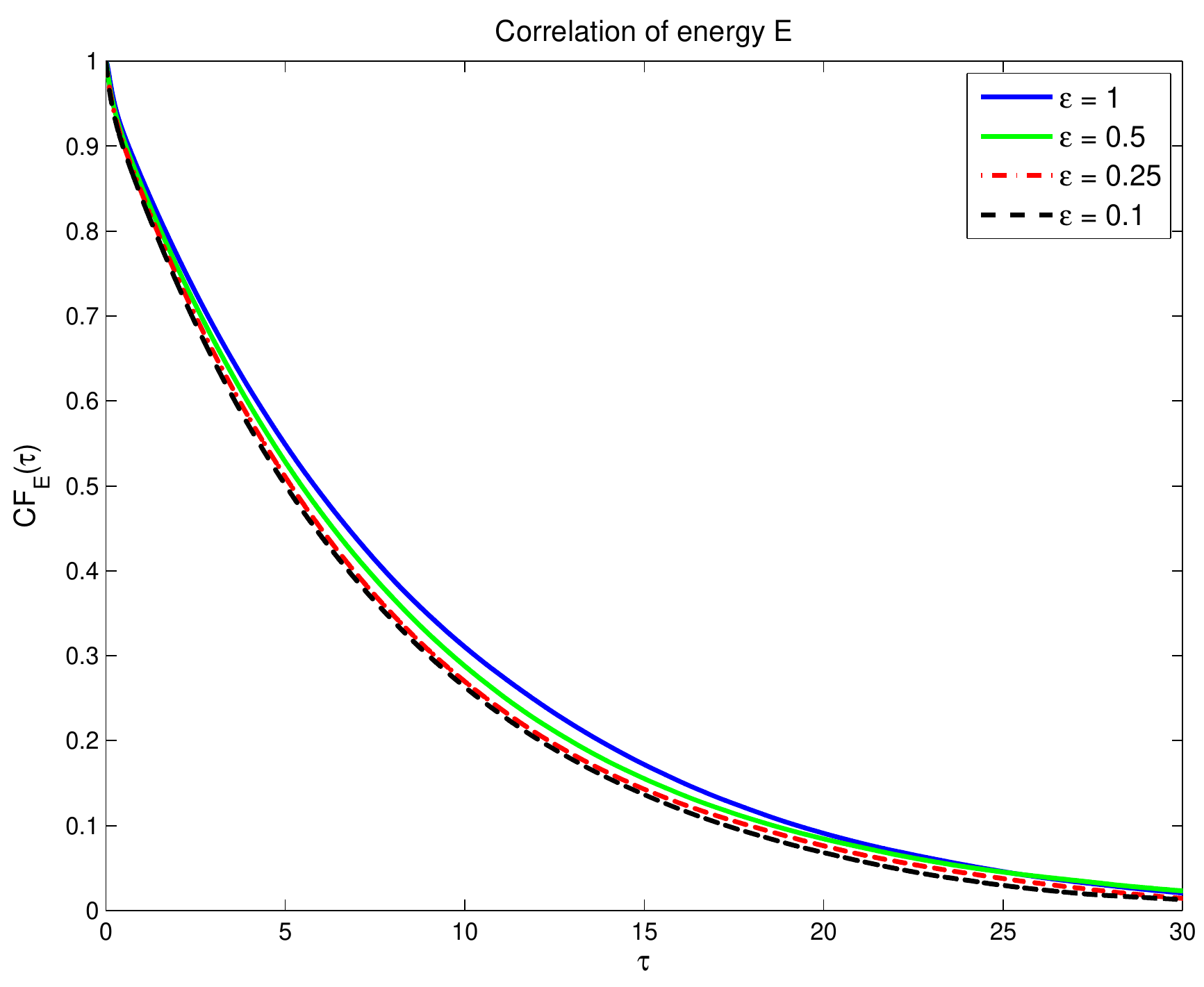}} 
\caption{Normalized correlation function of $E$ in that full model \eqref{triad} 
for $\eps=1, \, 0.5, \, 0.25, \, 0.1$.}
\label{fig:cfefull}
\end{figure}

Next, we examine the behavior of the full model as $\eps\to0$ numerically. In particular, 
we consider four values of $\eps$
\[
\eps = 1, \, 0.5, \, 0.25, \, 0.1
\]
and perform simulations to illustrate changes in the statistical behavior of the full model as $\eps\to 0$.
The stationary density for $x$ and $E$ agrees very well with the analytical predictions in \eqref{rhox}, \eqref{rhoE}
for any $\eps$ and we only depict the stationary correlation functions for $x$ and $E$ in Figures \ref{fig:cfxfull}
and \ref{fig:cfefull}, respectively. The stationary correlation function for a stochastic variable $z(t)$ is given by
\[
CF_z(\tau) = \bE z(t)z(t+\tau)  - (\bE z)^2
\]
and we normalize correlation functions by the variance, so that $CF(0)=1$.
Numerical simulations show that correlation functions of both, $x$ and $E$, 
change only slightly as $\eps\to 0$. Therefore, we expect that the stochastic mode reduction will 
be in a good agreement with the full model even for $\eps=1$.

\begin{figure}[ht]
\centerline{
\includegraphics[width=12cm]{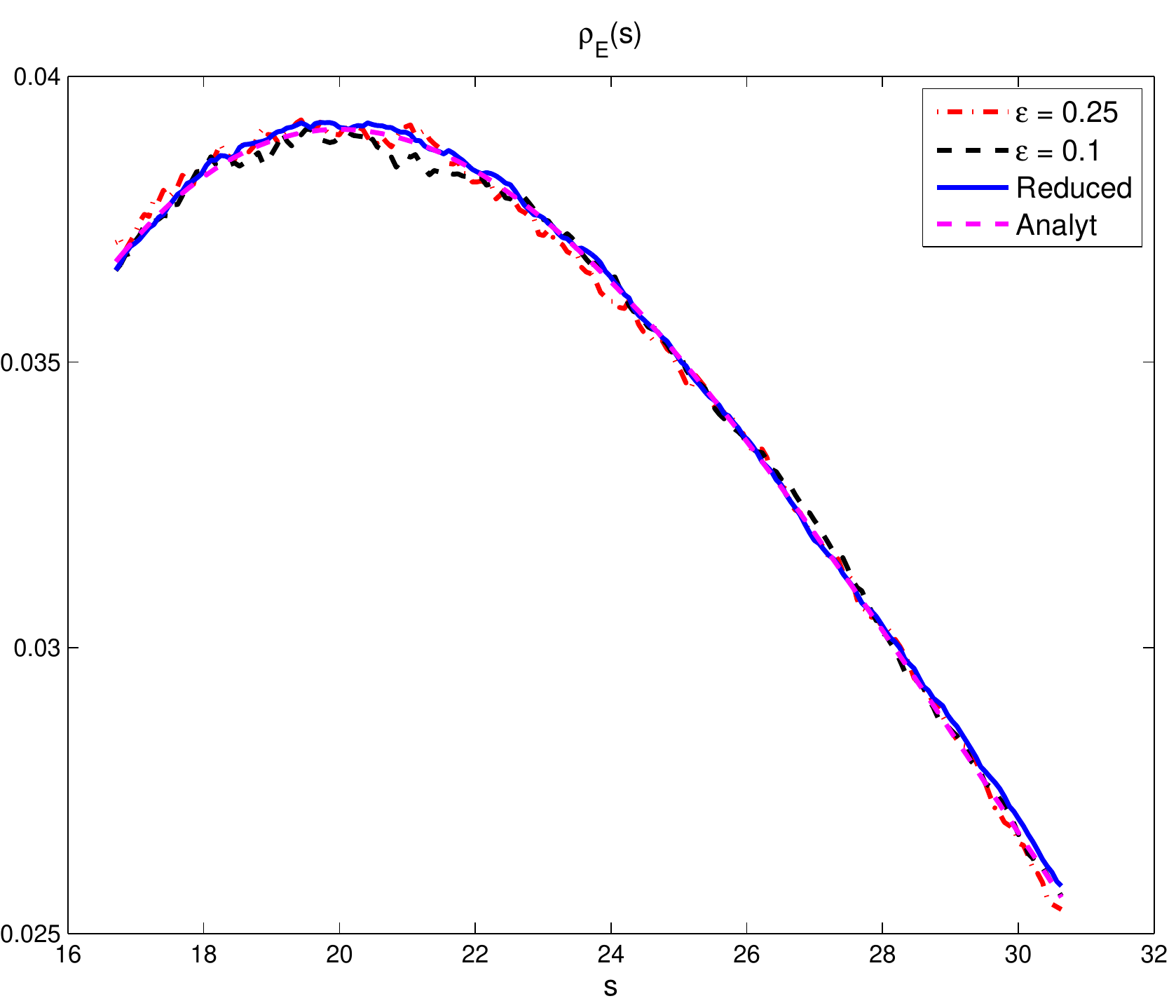}
}
\caption{Marginal Density of the slow variable $E$ in the full model with
$\eps=0.25$, $0.1$, Reduced Model, and the analytical formula in \eqref{rhoE}.}
\label{fig:pdf}
\end{figure}

As a final step, 
we compare the equilibrium statistical behavior of the full and reduced models to validate the 
performance of the stochastic mode reduction.
As expected, the marginal densities of $x$ and $E$ agree quite well in all simulations.
Here we present only the numerical results for the marginal density of $E$, the marginal densities for $x$
computed numerically completely overlap with each other and agree perfectly with the analytical 
prediction in \eqref{rhox}.
Comparison of the marginal density  $E$ is presented in Figure \ref{fig:pdf}.

Comparison of the two-point statistical quantities is a much more severe test for the mode-reduction 
since the mode reduction does not, in general, preserves the time-scales of the slow variables and the two-point
statistics for the full and reduced models can fail to agree due to the insufficient scale separation or
invalidity of the underlying assumptions (e.g. ergodicity, uniform distribution on the sphere, etc.).

\begin{figure}[ht]
\centerline{
\includegraphics[width=8cm]{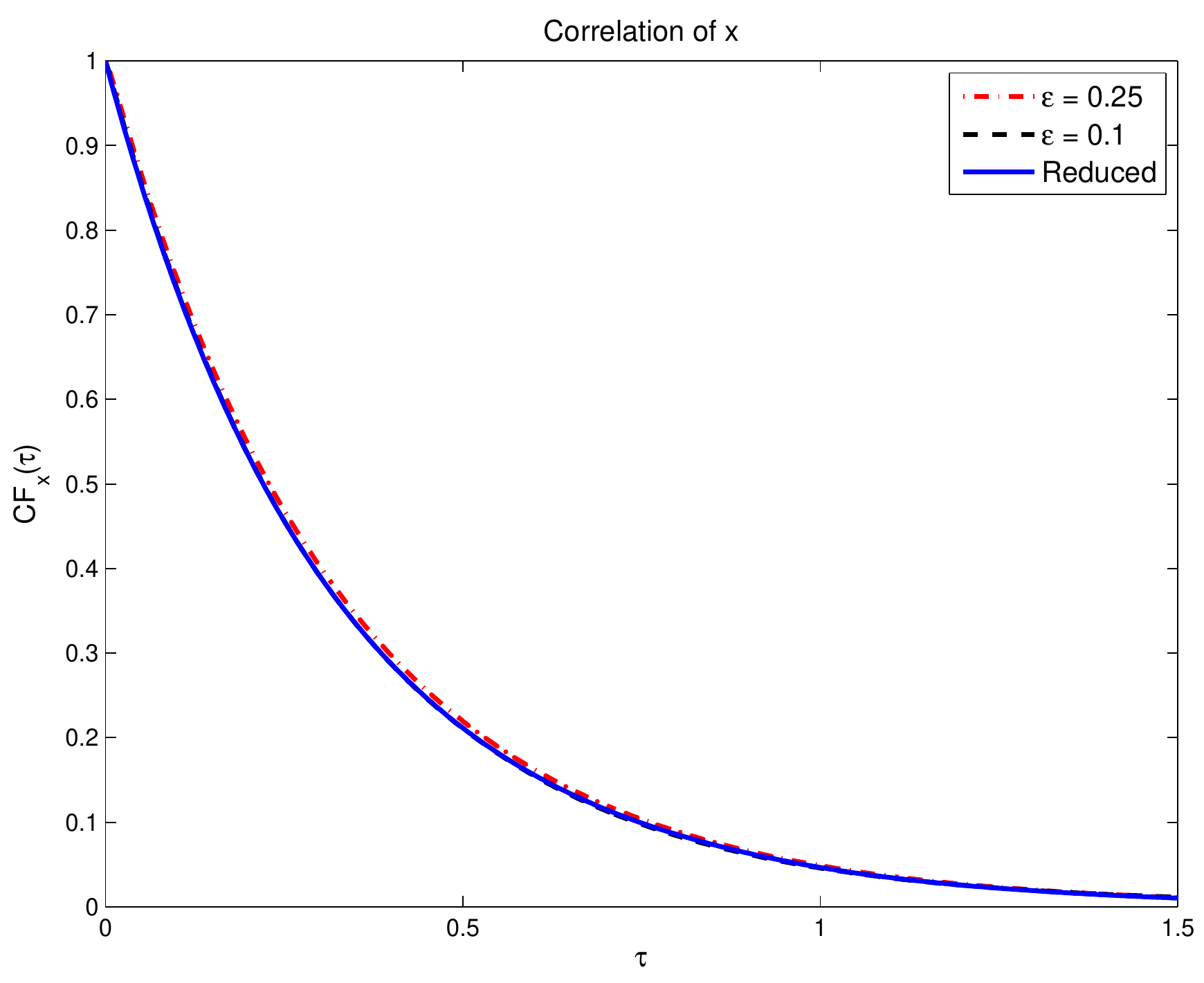}
\includegraphics[width=8cm]{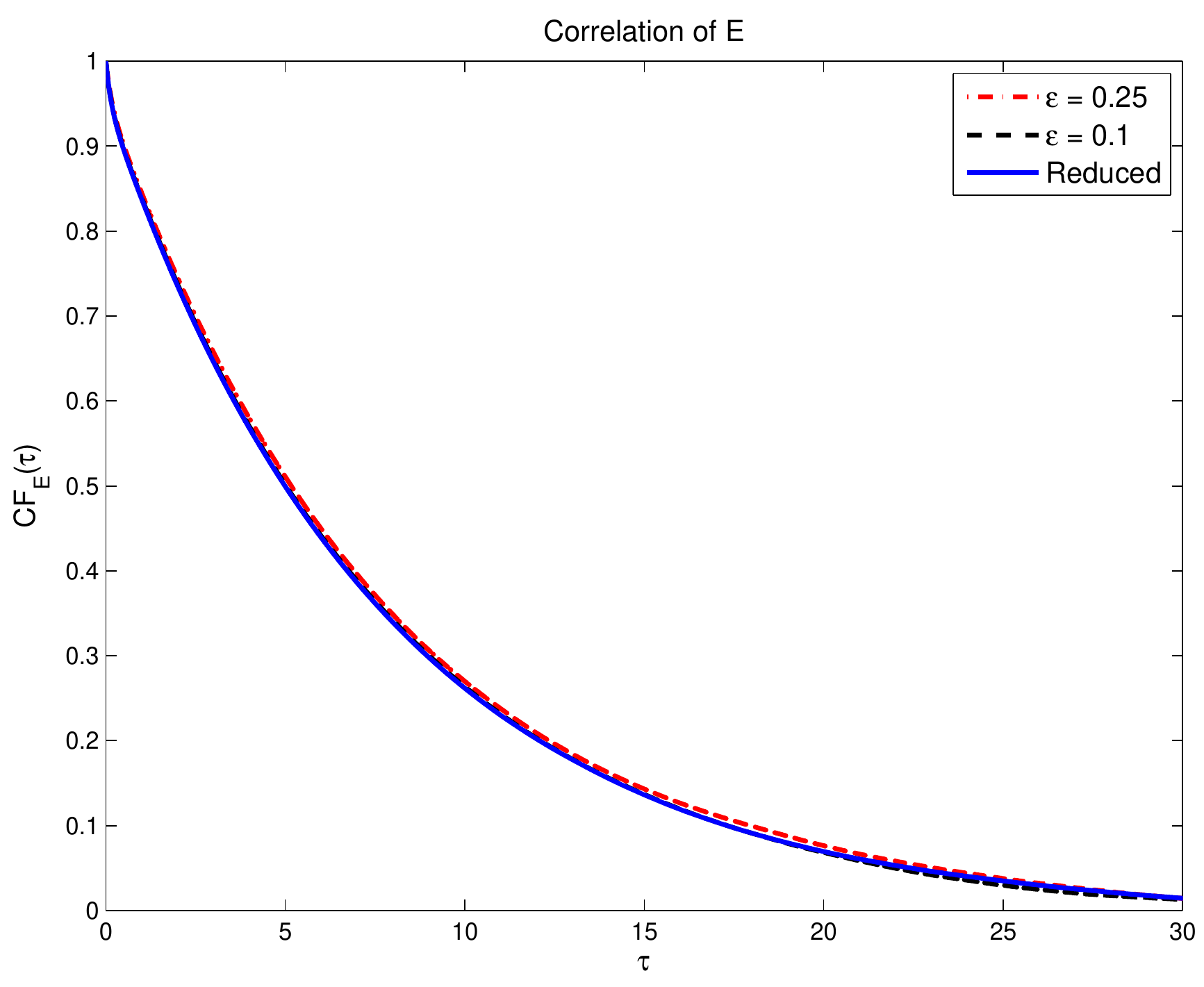}}
\caption{Correlation function of $CF_x$ and $CF_E$ in the full model \eqref{triad} with $\eps=0.25$, $0.1$ 
and in the Reduced model \eqref{reduced}. The solid blue line (reduced model) overlaps completely 
with the black dashed line (full mode with $\eps=0.1$),}
\label{fig:cf}
\end{figure}

Since the two-point statistical quantities cannot be computed analytically, 
we compare the results of numerical simulations of the full model for $\eps=0.25, \, 0.1$ and the reduced model.
Comparison of the correlation functions for $x$ and $E$ is presented in Figure \ref{fig:cf}.
Correlation functions for both, $x$ and $E$ in the simulations of the full model agree very well with the 
numerical results from the reduced model. This demonstrates the validity of the mode reduction procedure 
in the setup considered in this paper, i.e. when it is necessary to consider the energy of the fast sub-system
as an additional slow variable.

\begin{figure}[ht]
\centerline{
\includegraphics[width=8cm]{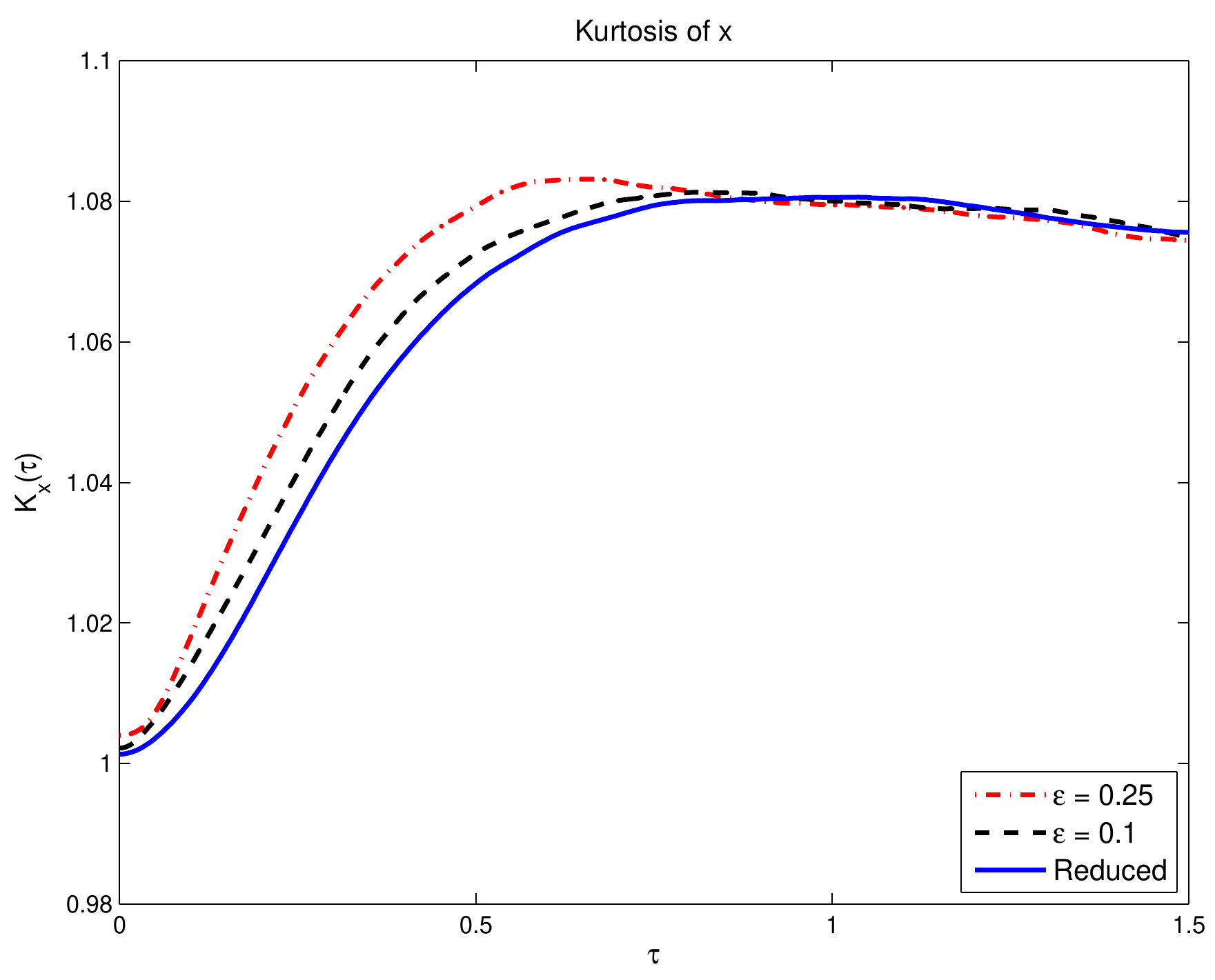}
\includegraphics[width=8cm]{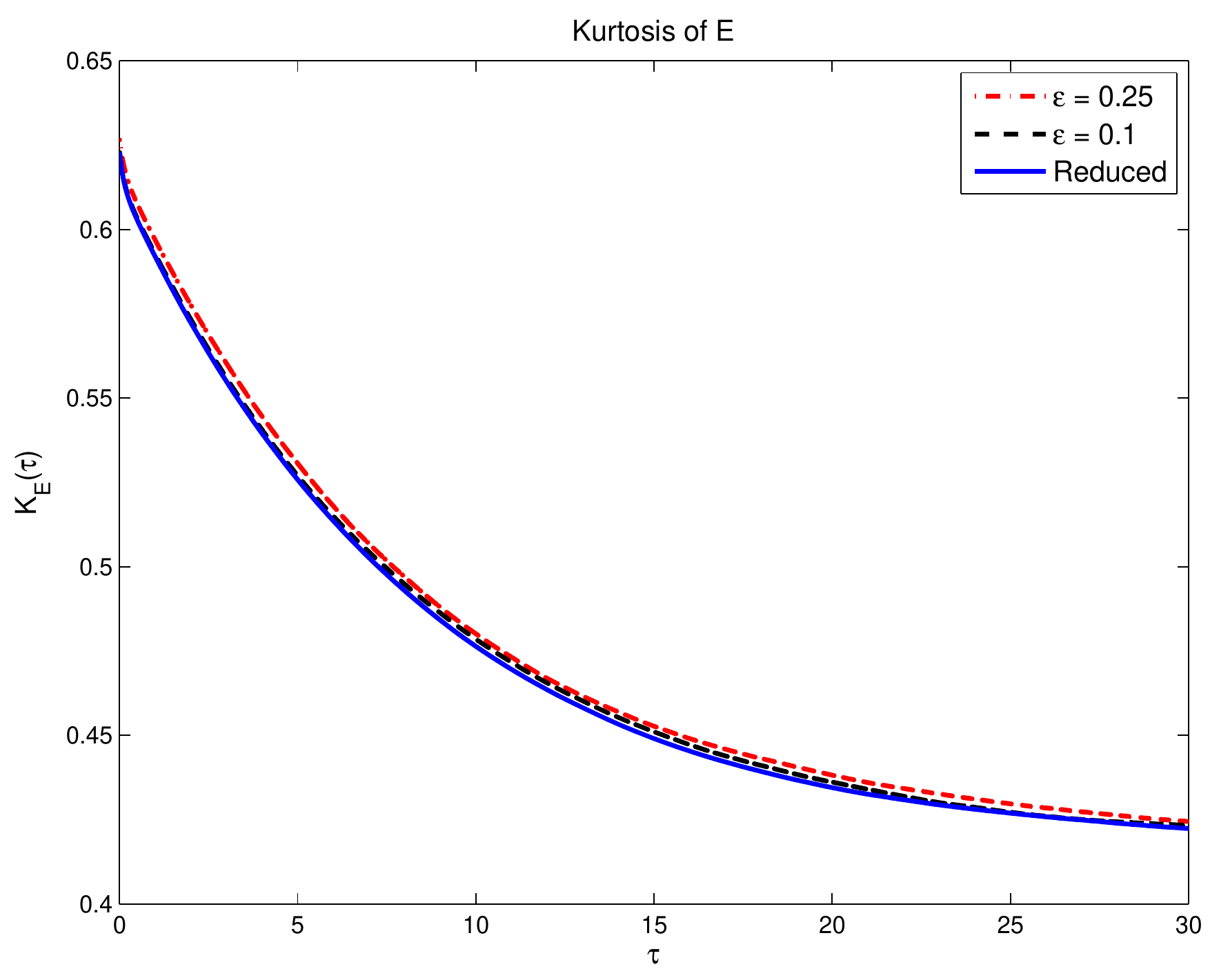}
}
\caption{Kurtosis \eqref{kurt} of slow variables $x$ and $E$ in full model with $\eps=0.25$, $0.1$ and in the Reduced model. The vertical scale for the Kurtosis of $x$ is rather small; errors between the full and reduced model for the Kurtosis of $x$ are approximately 1\%.}
\label{fig:kurt}
\end{figure}

It is also important to compare higher-order statistical quantities to analyze the non-Gaussian behavior of the system. 
Since the equilibrium distribution of $x$ is Gaussian and the correlation function
exhibits exponential decay,  one can be easily deceived to believe that a linear Ornstein-Uhlenbeck model would be sufficient to
accurately reproduce the statistical behavior of $x$. To illustrate the importance of the functional form of the reduced model in 
\eqref{reduced}, we compare the behavior of the fourth-order two-point moment
\begin{equation}
\label{kurt}
K_z(\tau) = \frac{\bE z^2(t+\tau) z^2(t) }{\left( \bE z^2 \right)^2 + 2 \left( \bE z(t)z(t+\tau) \right)^2}
\end{equation}
in Figure \ref{fig:kurt}. For Gaussian processes $K_z(\tau) \equiv 1$ for all lags $\tau$. 
Therefore, $K_x(\tau)$ and $K_E(\tau)$ measure the non-Gaussian features of the corresponding stochastic processes.
$K_x(\tau)$ in the left part of Figure \ref{fig:kurt} demonstrates the non-Gaussian feature of the stochastic variable $x(t)$.
Moreover, this Figure also suggests that the scale-separation affects the convergence of the higher-order statistics; while correlation functions of $x$ in the full model for $\eps=0.25,\, 0.1$ and the reduced model 
nearly overlap (Figure \ref{fig:cf}), there is about 1\% discrepancy between the full model with $\eps=0.25$ and $\eps=0.1$ and, also, between the reduced model and the full model with $\eps=0.1$. The discrepancy with the full model with $\eps=1$ and $\eps=0.5$ (not shown here) is even more severe. On the other hand, the kurtosis for the slow variable $E$ in the reduced stochastic model \eqref{reduced} is in a much better agreement with the full model \eqref{fullmodel}  for all values of $\eps$. 

Although the energy variable $E$ is much slower than $x$ in the reduced model \eqref{reduced}
(e.g. compare the correlation functions $CF_E$ and $CF_x$ depicted in Figure \ref{fig:cf})
for the particular choice of interaction coefficients discussed in this example,
the structure of the reduced equation \eqref{reduced} does not allow to 
average the behavior of $E$ with respect to the stationary measure of $x$.
In particular, the structure of the diffusion term precludes one from further accelerating the 
$x$ variable without affecting the joint stationary distribution for $x,E$.
On the other hand, by the same argument, it would be impossible to average the behavior of
$x$ with respect to the stationary measure for $E$, if $E$ was the faster variable.
Therefore, the system in \eqref{reduced} is the correct reduced model regardless of the 
relative behavior of $x$ and $E$.


\section{Conclusions}

We have presented an extension of the stochastic mode-reduction
  strategy for systems with stochastically perturbed slow variables
  and energy-conserving fast sub-system. The coupling between the slow
  and fast variables induces a slow energy exchange between the slow
  variables and fast sub-system, and it is necessary to include the
  slowly-evolving energy of the fast sub-system as an additional
  variable in the reduced description of the full model. Our numerical
  simulations demonstrate that the proposed formalism captures the
  equilibrium statistical behavior of the slow variables very well.
  Moreover, our example also demonstrates that the ``hidden'' energy
  variable can evolve on a much slower time-scale than the slow
  variables in the full model and, thus, the statistical behavior of
  these slow variables cannot be averaged over the distribution of the
  energy.

  The derivation of the reduced equation follows a standard formalism
  which relies on the asymptotic expansion of the backward equation,
  and the reduced equation involves constants which describe some bulk
  equilibrium statistical properties of the fast sub-system. Although
  these constants cannot be deduced analytically, they can be computed
  from a single microcanonical simulation of the fast sub-system on a
  shell of constant energy.  Such simulation is not particularly
  demanding since it is performed on the natural time scale of the
  fast sub-system, which is not multiscale.

  The formalism presented in this paper extends the applicability of
  the stochastic mode-reduction strategy.  Moreover, this formalism
  can be applied to purely deterministic systems with chaotic behavior
  of the slow variables.  Such systems can arise, for instance, as
  spectral truncations of coupled hyperbolic-parabolic PDE models,
  which will be considered in future work.

\section{Acknowledgement}
The research if I.T. was partially supported by the NSF Grant
DMS-1109582. I.T. also acknowledges support from SAMSI and IMA as a
long-term visitor in the Fall of 2011 and Fall 2012, respectively. The
research of E. V.-E. was partially supported by the DOE ASCR
Grant DE-FG02-88ER25053, the NSF Grant DMS07-08140, and the ONR Grant
N00014-11-1-0345.

\clearpage

\appendix

\section{Derivation of the Reduced Operator}
\label{a1}

Here, we outline the derivation of the reduced backward operator in \eqref{L final}.
Substituting $L_1$ into \eqref{L} and neglecting derivatives $\partial_{y_j}$ to the right of $L_2^{-1}$
(since the effective operator is applied to the function $u_0(x,E,t)$
which does not involve the $y$-variables), we obtain the following expression for $L$ 
\begin{eqnarray} 
& & - \bP L_1 L_2^{-1} L_1 = 
\label{Leff} \\ 
&& - \sum\limits_{jkj'k'} \int (A^{xyy}_{1jk} y_j y_k \partial_x  + A^{yxy}_{j1k} x y_k \partial_{y_j}
 - 2 A^{xyy}_{1jk} x y_j y_k \partial_E )
L_2^{-1} 
\nonumber \\
&& \qquad \qquad 
(A^{xyy}_{1j'k'} y_{j'} y_{k'} \partial_x - 2 A^{xyy}_{1j'k'} x y_{j'} y_{k'} \partial_E ) \, d\mu(\vec{y}|E) =
\nonumber \\ 
&& -  \sum\limits_{jkj'k'} \int\limits_0^\infty \int (A^{xyy}_{1jk} y_j y_k \partial_x  + A^{yxy}_{j1k} x y_k \partial_{y_j}
 - 2 A^{xyy}_{1jk} x y_j y_k \partial_E )
\nonumber \\  
&& \qquad \qquad
(A^{xyy}_{1j'k'} y_{j'}(\tau) y_{k'}(\tau) \partial_x - 2 A^{xyy}_{1j'k'} x y_{j'}(\tau) y_{k'}(\tau) \partial_E )  \, d\mu(\vec{y}|E)  d\tau,
\nonumber
\end{eqnarray}
where we substituted the action of $L_2^{-1}$ which corresponds to propagating the initial conditions 
$y_j$ in time using fast sub-system and integrating over all times.
Here, $y_j(\tau)$ is the solution on the  fast sub-system on the energy level $E$.

\noindent
\textbf{Deriving Diffusion Terms.} \\
The diffusion part of the reduced operator can be computed by collecting second derivatives in \eqref{Leff}
\begin{eqnarray}
&& L_{diff} = 
\sum\limits_{jkj'k'} A^{xyy}_{1jk} A^{xyy}_{1j'k'}  \int\limits_0^\infty \int y_j y_k  y_{j'}(\tau) y_{k'}(\tau) \, \rho(\vec{y}|E) \, d\vec{y} d\tau
\nonumber \\
&& 
\qquad \qquad \qquad 
\left[ \partial^2_x - 2x \partial_x \partial_E - 2x \partial_E \partial_x + 4x^2 \partial^2_E \right],
\nonumber
\end{eqnarray}
where $\rho(\vec{y}|E)$ is given by \eqref{mu}.
Let the double-integral in the expression above be denoted by a constant $Q$. This constant describes some bulk statistical properties of the fast sub-system.
In particular, it is the area under fourth two-point moment in fast subsystem.
$$
Q = \sum\limits_{jkj'k'} A^{xyy}_{1jk} A^{xyy}_{1j'k'} 
\int\limits_0^\infty \int y_j y_k  y_{j'}(\tau) y_{k'}(\tau) 
\rho(\vec{y}|E) \, d\vec{y}  \, d\tau.
$$

\noindent
\textbf{Rescaling of the Fast subsytem.} \\
Note that fast subsystem \eqref{fastsubs} is a deterministic model which conserves the  $E$. 
On the other hand, the fast sub-system is invariant under the following transformation 
$$
y = \left(\frac{E}{n}\right)^{1/2} y^*, \quad t = \left(\frac{E}{n}\right)^{-1/2}t^*
$$
where $y^*$ evolves on the energy shell $E=n$.
Hence, the area under the fourth moment $Q$ can be rewritten as
\begin{eqnarray}
Q &=& \left(\frac{E}{n}\right)^{3/2} 
\sum\limits_{jkj'k'} A^{xyy}_{1jk} A^{xyy}_{1j'k'} 
\int\limits_0^\infty \int y_j^* y_k^*  y_{j'}^*(\tau) y_{k'}^*(\tau) 
\rho(\vec{y}|E=n) \, d\vec{y} \, d\tau
\nonumber \\
&=& \left(\frac{E}{n}\right)^{3/2} M,
\label{Q}
\end{eqnarray}
where $M$  is the area under fourth moment in the fast subsystem \eqref{fastsubs} with fixed energy $E = n$
\begin{equation}
M = \sum\limits_{jkj'k'} A^{xyy}_{1jk} A^{xyy}_{1j'k'} 
\int\limits_0^\infty \int y_j y_k  y_{j'}(\tau) y_{k'}(\tau) 
\rho(\vec{y}|E=n) \, d\vec{y} \, d\tau .
\label{M1}
\end{equation}
The constant $M$ cannot be derived analytically. Therefore, we have to estimate the constant $M$ numerically.
However, the constant $M$ has to be computed only once from the numerical simulations of the fast sub-system alone,
which is not a multiscale model.

\noindent
\textbf{The Diffusion Matrix.} \\
Let $\Sigma$ be the diffusion matrix in the backward equation.. 
The square of the diffusion matrix for reduced model can be written as
 \[
\Sigma^2 = 
2 Q \left( \begin{array}{cc}
1 & -2x \\
-2x & 4x^2 \\
\end{array} \right) = 
2 M \left(\frac{E}{n}\right)^{3/2} \left( \begin{array}{cc}
1 & -2x \\
-2x & 4x^2 \\
\end{array} \right).
\]

\noindent
\textbf{Deriving the Drift and Diffusion terms using uniform measure on the sphere.}\\
The derivation of the drift terms in the reduced operator can proceed in two possible ways.
First, since the diffusion in the reduced operator is known and the stochastic mode reduction preserves the
invariant distribution of the slow variables \eqref{rhox}, \eqref{rhoE}, one can derive the drift terms
in the Fokker-Planck operator
which are consistent with the stationary distribution of $x$ and $E$ in the full model.
This derivation relies heavily on the known form of the stationary distribution for the
slow variables and, thus, the particular form of the self-interactions of the slow variables given by 
the operator $L_0$ and the form of the coupling terms between the slow and fast variables.

As an alternative derivation, we proceed by formally manipulating the derivatives in the expression for the 
reduced operator \eqref{Leff} to obtain the drift terms. To this end, one also has to differentiate the
invariant measure $\mu(\vec{y}|E)$ in the expression \eqref{Leff}.

We assume that the fast subsystem is ergodic on the hypersphere defined by energy of the fast subsystem 
with respect to the uniform distribution on the sphere.
Then the stationary distribution of the fast subsystem can be written as \eqref{mu}
\begin{equation}
d\mu(\vec{y}|E) = S_n^{-1} E ^{1-n/2} \delta \left(E - E_0\right) \, d\vec{y},
\label{mu1}
\end{equation}
where $n$ is the total number of fast variables, $E_0$ is the given fixed energy in the fast subsystem, $S_n^{-1} E ^{1-n/2}$ is the normalizing factor such that $S_n$  does not depend on $E$, and $\delta()$ is the Dirac delta function.
The normalization factor $S_n^{-1} E ^{1-n/2}$ is derived from the condition
\[
\int d\mu(\vec{y}|E) = 1.
\]

Since the stationary measure of the fast subsystem depends on the energy $E$, the operator $L_2^{-1}$ also depends on the energy level $E$. Moreover,  $L_2^{-1}$  also depends on $y_i$ and the behavior of the time-dependent variables $y_j(\tau)$ depends on the energy level as well. 
Therefore, to understand the derivatives $\partial_E L_2^{-1}$ and $\partial_{y_i} L_2^{-1}$, we split the effective operator $L$ in \eqref{Leff} into three parts
\begin{eqnarray}
L &=& I_1 + I_2 + I_3,
\label{eff L splitting}
\end{eqnarray}
where 
\begin{eqnarray}
I_1 &=& -\int \rho(\vec{y}|E)  \sum_{j,k} A^{xyy}_{1jk} y_j y_k \partial_x  L_2^{-1} \left(\sum_{j,k} A^{xyy}_{1jk} y_j y_k \partial_x - 2 \sum_{j,k} A^{xyy}_{1jk} x y_j y_k \partial_E  
\right) d\vec{y}
\nonumber \\
&=&
\sum_{j,k} A^{xyy}_{1jk}
\int_0^{\infty} \int \rho(\vec{y}|E)   y_j y_k \partial_x \left(\sum_{j',k'} A^{xyy}_{1j'k'} y_{j'}(\tau) y_{k'}(\tau) [\partial_x - 2 x \partial_E] \right)  d\vec{y} d\tau ,
\label{I1}
\end{eqnarray}
\begin{eqnarray}
I_2 &=& -\int \rho(\vec{y}|E)  \sum_{j,k} A^{yxy}_{j1k} x y_k \partial_{y_j}  L_2^{-1} \left( \sum_{j,k} A^{xyy}_{1jk} y_j y_k \partial_x - 2 \sum_{j,k} A^{xyy}_{1jk} x y_j y_k \partial_E 
\right)  d\vec{y}
\nonumber \\
&=&
x \sum_{j,k,j',k'} A^{yxy}_{j1k} A^{xyy}_{1j'k'} \int _0^{\infty} \int \rho(\vec{y}|E)  y_k \partial_{y_j} \left(
y_{j'}(\tau) y_{k'}(\tau) \right) [\partial_x - 2 x \partial_E ] \, d\vec{y} d\tau,
\label{I2}
\end{eqnarray}
\begin{eqnarray}
I_3 &=& \int \rho(\vec{y}|E)  2 \sum_{j,k} A^{xyy}_{1jk} x y_j y_k \partial_E  L_2^{-1} \left( \sum_{j,k} A^{xyy}_{1jk} y_j y_k \partial_x - 2 \sum_{j,k} A^{xyy}_{1jk} x y_j y_k \partial_E 
\right) d\vec{y}
\nonumber \\
&=&
- 2 x \sum_{j,k,j',k'} A^{xyy}_{1jk} A^{xyy}_{1j'k'} \int_0^{\infty} \int \rho(\vec{y}|E) y_j y_k \partial_E \left( 
y_{j'}(\tau) y_{k'}(\tau) \, [\partial_x - 2 x \partial_E] \right) d\vec{y} d\tau.
\label{I3}
\end{eqnarray}
%

\noindent
\textbf{Computing $I_1$}\\
Since the fast sub-system does not depend on the value of the slow variable $x$, the calculation of $I_1$ is
straightforward and we obtain 
\begin{eqnarray} \nonumber 
I_1 &=& 
\sum_{j,k} A^{xyy}_{1jk} \sum_{j',k'} A^{xyy}_{1j'k'}
\int_0^{\infty} \int   y_j y_k  y_{j'}(\tau) y_{k'}(\tau)   \rho(\vec{y}|E)  \, d\vec{y} d\tau 
\partial_x [\partial_x - 2 x \partial_E] ,
\end{eqnarray}
where $y_j(\tau)$ is the solution of the fast subsystem with $y_j(0) = y_j$.
The expression for $I_1$ involves the area under the fourth two-points moment, similar to the discussion before
\[
I_1 = Q (\partial_{xx} - 2 x \partial_{xE} - 2 \partial_E ),
\]
where 
and $Q_{j,k,j',k'}$ is the integrated fourth-order two-point moment, in the fast subsystem \eqref{fastsubs} on the energy level $E$. $I_1$ can be further simplified after rescaling of the fast subsystem
to the energy shell $E=n$
\begin{eqnarray}
I_1 = \left(\frac{E}{n}\right)^{3/2} M \, \left[\partial_{xx} - 2 x \partial_{xE} - 2 \partial_E \right],
\label{I term}
\end{eqnarray}
where $M$ is given by \eqref{M}.

\noindent
\textbf{Computing $I_2$}\\
The operator $I_2$ involves differentiation with respect to $y_j$ which cannot be evaluated in a straightforward manner.
Instead, this differentiation can be switched onto the $\rho(\vec{y}|E)$ by integrating by parts.
Therefore, we obtain that $I_2$ can be expressed as
\begin{eqnarray}
I_2 &=& - 2 x S_n^{-1} E^{1-n/2}  J \, \left[\partial_x  - 2 x \partial_E \right],
\label{II term}
\end{eqnarray}
where $J$ is defined as follows
\begin{eqnarray}
J = \sum_{j,k,j',k'}  A^{yxy}_{j1k} A^{xyy}_{1j'k'} \int\limits_0^\infty \int y_j y_k  y_{j'}(\tau) y_{k'}(\tau) \delta' \left(E - \sum_{r=1}^n y_r^2 \right) \, d\vec{y} d\tau.
\label{Jdefn}
\end{eqnarray}
Note that the expression above involves derivative of the Dirac delta function.

\noindent
\textbf{Computing $I_3$}\\
We can use the product rule to express the $I_3$ as follows
\begin{eqnarray}\nonumber
I_3 &=&  - 2 x \sum_{j,k,j',k'} A^{xyy}_{1jk} A^{xyy}_{1j'k'} \int_0^{\infty} \int \rho(\vec{y}|E)  y_j y_k \partial_E \left[ y_{j'}(\tau) y_{k'}(\tau) \right] d\vec{y} d\tau ~ \left[\partial_{x} - 2 x \partial_{E}\right] 
\\ \nonumber
&&- 2 x Q \, \left[ \partial_{Ex} - 2 x \partial_{EE} \right],
\end{eqnarray}
where $Q$ is the integrated fourth-order two-point moment in \eqref{Q}.
Since the behavior of $y(\tau)$ depends on the energy level, $E$, the first term is non-zero.
Let the first term left in $I_3$ be denoted by $\tilde{I}_3$
\begin{equation}
\tilde{I}_3 = \int_0^{\infty} \int \mu(\vec{y}|E) y_j y_k \partial_E \left[ Y_{j'}(t) Y_{k'}(t)\right] \, d\vec{y} \, d\tau ,
\label{I 3 sub}
\end{equation}
then $I_3$ can be written as 
\begin{eqnarray}
I_3 &=&  - 2 x \sum_{j,k,j',k'} A^{xyy}_{1jk} A^{xyy}_{1j'k'} \tilde{I}_3 ~ [\partial_{x} - 2 x \partial_{E}] - 2 Q x  [\partial_{Ex} - 2 x \partial_{EE}].
\label{I 3 1}
\end{eqnarray}

In order to understand the structure if the term $\tilde{I}_3$ 
Let us compute the derivative of the area under the fourth-order moment.
\begin{eqnarray*}
&& \partial_E \int_0^{\infty} \int  y_j y_k y_{j'}(\tau) y_{k'}(\tau)  \rho(\vec{y}|E) \, d\vec{y} d\tau  =
\\ \nonumber
&& \int_0^{\infty} \int 
\biggl( y_j y_k y_{j'}(\tau) y_{k'}(\tau) \partial_E\left[\rho(\vec{y}|E)\right] +
       \rho(\vec{y}|E) y_j y_k \partial_E\left[ y_{j'}(\tau) y_{k'}(\tau)\right]  \biggr) \, d\vec{y} \, d \tau
\\ \nonumber
&=&
\int_0^{\infty} \int y_j y_k y_{j'}(\tau) y_{k'}(\tau) \partial_E\left[\rho(\vec{y}|E)\right] \, d\vec{y} \, d \tau + \tilde{I}_3 .
\end{eqnarray*}
Therefore, the required term can be rewritten as
\begin{eqnarray}
&& \sum_{j,k,j',k'} A^{xyy}_{1jk} A^{xyy}_{1j'k'}  \tilde{I}_3 = 
\\ \nonumber
&& \partial_E Q - 
\sum_{j,k,j',k'} A^{xyy}_{1jk} A^{xyy}_{1j'k'} \int_0^{\infty} \int  y_j y_k y_{j'}(\tau) y_{k'}(\tau) \partial_E \left[\rho(\vec{y}|E) \right] \, d\vec{y} d\tau .
\label{I 3 sub 1}
\end{eqnarray}
Since the stationary measure of the fast variables $\mu(\vec{y}|E)$ given by \eqref{mu}
we can compute the derivative of $\rho(\vec{y}|E)$ as follows
\begin{eqnarray*}
\partial_E \left [\rho(\vec{y}|E) \right] &=& (1-n/2) E^{-1} \rho(\vec{y}|E) + S_n^{-1} E^{1 - n/2} \delta' \left(E - \sum_{i=1}^n y_k^2\right),
\end{eqnarray*}
where $\delta'(\cdot)$ represent the distributional derivative of Dirac delta function.
Substituting the above expression into $\tilde{I}_3$ we obtain

\begin{eqnarray}
&& - 2 x \sum_{j,k,j',k'} A^{xyy}_{1jk} A^{xyy}_{1j'k'} \tilde{I}_3 =
\\ \nonumber
&& \qquad -2x \partial_E Q + x (2-n) E^{-1} Q + 2x S_n^{-1} E ^{1 - n/2} J,
\nonumber
\end{eqnarray}
where $Q$ and $J$ are given by \eqref{Q} and \eqref{Jdefn}, respectively. Therefore,
putting everything together, we obtain 
\begin{eqnarray}\nonumber
I_3 &=&  - x \frac{n+1}{n^{3/2}} E^{1/2} M  (\partial_{x} - 2 x \partial_{E})  - 2 x \frac{E^{3/2}}{n^{3/2}} M (\partial_{Ex} - 2 x \partial_{EE}) \\ &&
+ 2 x S_n^{-1} E ^{1 - n/2} J (\partial_{x} - 2 x \partial_{E}).
\label{III term}
\end{eqnarray}

Substituting $I_1$, $I_2$, and $I_3$ into the effective operator we obtain that the terms involving $J$ (derivative $\delta'()$) cancel, and the effective operator becomes
\begin{eqnarray} 
L= &=& - \frac{E^{1/2}}{n^{3/2}} M  \left( (n+1) x \partial_{x} + 2 E \partial_E - 2 x^2 (n+1) \partial_{E} \right) 
\nonumber \\ 
& & + \frac{E^{3/2}}{n^{3/2}} M \left(\partial_{xx} - 2 x \partial_{xE} - 2 x \partial_{Ex} + 4 x^2 \partial_{EE} \right),
\label{L final1}
\end{eqnarray}
where $M$ is summations of the area under fourth-order two-point moment given by \eqref{M}. 
The diffusion part of the operator agrees with the previously computed matrix, i.e.
the diffusion part of $L$ can be rewritten as
\begin{equation}
L_{diff} = \frac{1}{2} \sum_{i,k=1,2} (D D^T)_{ik} \partial_{z_i} \partial_{z_k}
\label{L diff} 
\end{equation}
where $(z_1, z_2) \equiv (x,E)$ and 
\begin{equation}
\label{D}
D = \sqrt{2M} \left(\frac{E}{n}\right)^{3/4}
\left( \begin{array}{cc}
1 & 0\\
-2x & 0\\
\end{array} \right).
\end{equation}

Using the effective $L$ operator in \eqref{L final1} and diffusion in \eqref{D}, we obtain a reduced SDE model for $x$ 
and $E$ given by
\begin{eqnarray}
dx &=& - \gamma x dt - (n+1) M x \frac{E^{1/2}}{n^{3/2}} + \sigma dW_1 + 
\sqrt{2 M} \left(\frac{E}{n}\right)^{3/4} dW_2 ,
\nonumber \\
dE &=&  -2 M \left(\frac{E}{n}\right)^{3/2} + 2 (n+1) M x^2 \frac{E^{1/2}}{n^{3/2}} - 2 \sqrt{2 M}x \left(\frac{E}{n}\right)^{3/4} dW_2 ,
\nonumber
\end{eqnarray}
where $M$ is a numerical constant \eqref{M} which needs to be computed numerically; $M$ represents integrated fourth-order two-point moment computed in the fast sub-system \eqref{fastsubs} on the energy level $E=n$.

\clearpage

\section{Interaction Coefficients in the Full Model}
\label{a2}

\begin{table}[ht]
	\centering
\begin{tabular}{|c|c|c|c|}
	\hline
	$(j,k)$          & $A_{1jk}$ & $A_{j1k}$ & $A_{k1j}$ \\ \hline
    ${1,2}$           &   $1.2$ & $-0.55$ & $-0.65$\\ \hline  
    ${8,9}$       &   $0.525$ & $0.25$ & $-0.775$\\ \hline
    ${4,10}$         &   $1.35$ & $-0.725$ & $-0.625$\\ \hline
    ${5,6}$           &   $1.125$ & $-0.5$ & $-0.625$\\ \hline
    ${3,7}$           &   $1.35$ & $-0.725$ & $-0.625$\\ \hline
    ${1,10}$         &   $0.525$ & $0.25$ & $-0.775$\\ \hline
    ${2,4}$           &   $1.2$ & $-0.55$ & $-0.65$\\ \hline
    ${5,8}$           &   $1.125$ & $-0.5$ & $-0.625$\\ \hline
    ${7,9}$           &   $0.875$ & $-0.3$ & $-0.575$\\ \hline
    ${3,6}$           &   $1.25$ & $-0.625$ & $-0.625$\\ \hline
\end{tabular}
\caption{Coefficients $A^{xyy}_{1jk}, A^{yxy}_{j1k}, A^{yyx}_{jk1}$ used in coupling of $x$ and $y$ variables in the full model \eqref{triad}.}
\label{tabA}
\end{table}

\begin{table}[ht]
	\centering
\begin{tabular}{|c|c|c|c|}
	\hline
	$(i,j,k)$        &   $B_{ijk}$ & $B_{jki}$ & $B_{kij}$ \\ \hline
    ${1,2,3}$        &   $2$ & $2.5$ & $-4.5$\\ \hline
    ${1,2,4}$        &   $4.2426$ & $2.8284$ & $-7.071$\\ \hline
    ${1,2,9}$        &   $-1.2247$ & $2.9393$ & $-1.7146$\\ \hline 
    ${1,2,10}$      &   $2.1166$ & $2.9103$ & $-5.0269$\\ \hline
    ${1,3,4}$        &   $1.7321$ & $2.5981$ & $-4.3302$\\ \hline
    ${1,5,6}$        &   $3.8013$ & $4.9193$ & $-8.7206$\\ \hline
    ${1,9,10}$      &   $3.9598$ & $-2.2627$ & $-1.6971$\\ \hline 
    ${2,3,4}$        &   $-2$ & $4$ & $-2$\\ \hline
    ${2,5,6}$        &   $-4.5$ & $2.1$ & $2.4$\\ \hline
    ${2,9,10}$      &   $1.7393$ & $1.4230$ & $-3.1623$\\ \hline
    ${3,7,8}$        &   $1.1608$ & $2.3217$ & $-3.4825$\\ \hline
    ${4,7,8}$        &   $-1.7321$ & $-2.0785$ & $3.8106$\\ \hline
    ${5,6,7}$        &   $2.9566$ & $2.0912$ & $-5.0478$\\ \hline
    ${5,6,8}$        &   $-2.6192$ & $-1.4966$ & $4.1158$\\ \hline
    ${5,7,8}$        &   $4.6476$ & $2.7111$ & $-7.3587$\\ \hline 
    ${5,6,9}$        &   $-3$ & $-1.8$ & $4.8$\\ \hline
    ${5,6,10}$      &   $1.8554$ & $2.2677$ & $-4.1231$\\ \hline
    ${6,7,8}$        &   $4.6669$ & $2.9698$ & $-7.6367$\\ \hline
    ${8,9,10}$      &   $3.923$  & $2.3974$ & $-6.3204$\\ \hline
\end{tabular}
\caption{Coefficients $B^{yyy}_{ijk}$ used in coupling of $y$ variables in the
full model \eqref{triad}.} 
\label{tabB}
\end{table}


\end{document}